\documentclass[12pt]{amsart}
\usepackage[centertags]{amsmath}
\usepackage{graphicx}
\usepackage{amsfonts,amssymb,latexsym,epsfig,amsthm}
\usepackage{fullpage}

\newtheorem{theorem}{Theorem}

\newtheorem{claim}[theorem]{Claim}

\newtheorem{conjecture}[theorem]{Conjecture}

\newtheorem{definition}[theorem]{Definition}

\newtheorem{lemma}[theorem]{Lemma}

\newtheorem{proposition}[theorem]{Proposition}
\newtheorem{remark}[theorem]{Remark}

\input{tcilatex}

\begin{document}
\title{On a multi-particle Moser-Trudinger Inequality}
\author{Hao Fang}
\address{Courant Institute of Mathematical Sciences, New York University}
\date{2003}
\maketitle

\begin{abstract}
We verify a conjecture of Gillet-Soul\'{e}. We prove that the determinant of
the Laplacian on a line bundle over $\mathbb{CP}^{1}$ is always bounded from
above. This can also be viewed as a multi-particle generalization of the
Moser-Trudinger Inequality. Furthermore, we conjecture that this functional
achieves its maximum at the canonical metric.\ We give some evidence for
this conjecture, as well as links to other fields of analysis.
\end{abstract}

\section{Introduction}

\setcounter{equation}{0} \setcounter{theorem}{0}For a compact complex curve $%
C$ with a Hermitian metric, $g^{C},$ and a line bundle $l$ over $C,$ with
associated Hermitian metric, $g^{l}$, the determinant of the $\overline{%
\partial }-$Laplacian over $\mathcal{O}(l)$ is defined by the method of
zeta-regularization. We denote its logarithm as $A(g^{C},g^{l}),$ a
functional of $g^{C}$ and $g^{l}.$

In~\cite{GiSo1},\ Gillet and Soul\'{e} posed the following:

\begin{conjecture}[Gillet-Soul\'{e}]
\label{conjecture} $A(g^{C},g^{l})$ is bounded from above by a constant
independent of choice of $g^{C}$ and $g^{l}$.
\end{conjecture}

Conjecture~\ref{conjecture} is motivated by the Arithmetic Riemann-Roch
theorem; its bounded-ness will imply the one-side bounded-ness of the
modified arithmetic Betti number, first introduced by Gillet and Soul\'{e}.

From an analytic point of view, $A(g^{C},g^{l})$ is a spectral invariant
that is naturally linked to the\ quantization problem of Toeplitz
operators.\ The conjecture can alo be viewed a natural holomorphic extension
of the classical\ Szeg\"{o} limit theorem, which was originally stated for $%
S^{1}$~\cite{Szego-1}. See Section 5 for more details.

On each smooth closed curve $C,$ there exists a canonical metric with
constant curvature, $g_{0}^{C},$ which in turn induces a canonical metric on
$l,$ $g_{0}^{l}.$ By the Uniformization Theorem, general metrics on $C$ and $%
l$ are conformal to $g_{0}^{C}$ and $g_{0}^{l}$, respectively. According to
the Serre Duality, one only needs to consider ample line bundles over the
curve; hence, it also can be intepreted as an analytic torsion assoicated to
the line bundle $l\ $(See Section 2 for more details). A simple application
of the family Riemann-Roch theorem of Bismut-Gillet-Soul\'{e} \cite{BGS1}
then gives the anomaly formula of the regularized determinant in terms of
the anomaly of the Quillen metric and the L$^{2}$-metric. Furthermore, by a
similar computation, $g^{C}$ can be fixed (see \cite{GiSo1} for details);
hence, it is enough to prove the conjecture for $g^{C}=g_{0}^{C}$ and $%
g^{l}=\exp \varphi $ $g_{0}^{l}.$

In the same paper, Gillet and Soul\'{e} considered the case of $C$ being
rational.\ They proved the conjecture under the condition that the metric $%
g^{l}$ is endowed with an additional rotational symmetry. The proof is a
delicate extension of the original proof of Moser for the famous
Moser-Trudinger Inequality for two-sphere~\cite{Moser-1}.

In this paper, we remove the special condition and prove the conjecture for
the rational curve case in general:

\begin{theorem}
\label{main}Conjecture \ref{conjecture} holds for any holomorphic line
bundle over $\mathbb{CP}^{1}.$
\end{theorem}

Naturally, we consider the problem of exact upper bound. We conjecture that

\begin{conjecture}
\label{conjecture-2.0}$A(g_{0}^{C},g^{l})$ achieves its sharp upper bound
only when $g^{l}$ is the standard metric $g_{0}^{l}.$
\end{conjecture}

See Conjecture \ref{conjecture-2} for an equivalent but more concrete
statement.

By Theorem \ref{main-2}, Conjecture \ref{conjecture-2.0} holds for function $%
\varphi $ with large L$^{2}$ norm. We also show that constant functions are
local maximal points of the determinant functional $A;\mathit{\ }$\textit{%
i.e.}, Conjecture \ref{conjecture-2.0} holds for function $\varphi $ with
very small energy.

In this paper, we will give several interpretations of this conjecture.

Firstly, it can be viewed as sharp form of a multi-particle generalization
of the Moser-Trudinger inequality~\cite{Moser-1}. In fact, Theorem \ref{main}
can be viewed as the first step of an attempt to emulate Moser's original
proof of the original Moser-Trudinger Inequality; the special case solved
by\ Gillet and Soul\'{e} in \cite{GiSo1} is a delicate extension of the
Moser's original approach. We remark that Jost-Wang considered a different type of extension of
the Moser-Trudinger inequality relating to the Toda system~\cite{Jost-Wang}.

Secondly, it is a $S^{2}$ holomorphic extension of the classical Szeg\"{o}
limit theorem, which was stated for $S^{1}$ in the frame work of Fourier
analysis$.$ Notice that Okikiolu extended a weak form of Szeg\"{o} limit
theorem to $S^{2}$ and $S^{3}~$\cite{ko}. However, our conjecture is a sharp
inequality for all finite natural number $n$, and we consider a holomorphic
version. See Section 5 for more details.

Thirdly, it can be viewed as a problem related to the quantization of the
Toeplitz operator on $S^{2}.$ It is thus closedly related to earlier works
of Boutet de Monvel-Guillemin~\cite{BdM-G}\ and Uribe~\cite{uribe-1}.

Finally, it is a limit theorem when interpreted as a configuration problem
of random variables on $S^{2}.$\ This is partly inspired by the Kac's
probalistic approach to the Szeg\"{o} limit theorem~\cite{Kac-2}.

This paper is organized as follows. In Section 2, we set up the problem and
state the result of\ Gillet-Soul\'{e}. In Section 3, we prove Theorem~\ref%
{main}. In Section 4, we establish the local maximality of the constant
function for functional $A.$ In Section 5, we pose the new conjecture on the
sharp bound and give several formulation of the problem, relating it to
different classical results.

\noindent \textbf{Acknowledgement} The author would like to thank Paul Yang
for bringing this problem to his attention.\ Thanks are also due to Sun-Yung
Chang, Xianzhe Dai, Kate\ Okikiolu and Paul Yang for discussions. Part of
this work was done during the author's stay in MSRI in 2001.

\section{Set-up\setcounter{equation}{0} \setcounter{theorem}{0}}

In this Section, we define the related geometric objects and state the
previous results of \cite{GiSo1}.

From now on, we fix $C=\mathbb{CP}^{1}=\mathbb{C}\cup \{\infty \}$ with $z$
being the complex coordinate function. Let%
\begin{equation}
\mu ={\frac{\sqrt{-1}{dzd}\overline{z}}{{2\pi }\left( 1+\left\vert
z\right\vert ^{2}\right) {^{2}}}}
\end{equation}%
is the K\"{a}hler form of the standard Fubini-Study metric such that $%
\int_{C}\mu =1$.

Notice that, with our notation, the celebrated Moser-Trudinger Inequality
reads as~\cite{Moser-1}%
\begin{equation}
\log \int_{C}\exp (\varphi )\mu -\frac{1}{4}\int_{C}|\nabla \varphi |^{2}\mu
-\int_{C}\varphi \mu \leq 0.  \label{Moser-Trudinger}
\end{equation}

Fix a positive integer $n$. $\mathcal{O}(n)$ denotes the holomorphic line
bundle over $C$ of degree $n$. For a real $C^{\infty }$ function $\varphi $
over $C$, we can define the following metric $g_{\varphi }$ on $\mathcal{O}%
(n)$:
\begin{equation}
<s,t>_{\varphi }=\int_{C}\frac{s\overline{t}}{{(1+}\left\vert z\right\vert {%
^{2})^{n}}}{(\exp {\varphi })\mu ,}
\end{equation}%
where $s$ and $t$ are two sections of $\mathcal{O}(n)$. We then have the
relation $g_{\varphi }=(\exp \varphi )g_{0}.$

We consider the following chain complex
\begin{equation}
0\rightarrow \Lambda ^{0,0}(C,\mathcal{O}(n))\overset{\overline{\partial }}{%
\rightarrow }\Lambda ^{0,1}(C,\mathcal{O}(n))\rightarrow 0.  \label{chain}
\end{equation}

The metric $g_{\varphi }$, and the constant curvature metric of $C,$ $%
g_{0}^{C}$, naturally induces metrics on $\Lambda ^{0,0}(C,\mathcal{O}(N))$
and $\Lambda ^{0,1}(C,\mathcal{O}(N))$. Hence, we define the adjoint of $%
\bar{\partial},$ $\bar{\partial}_{\varphi }^{\ast }$. Now we define
\begin{equation}
\Delta _{\varphi }=\bar{\partial}\bar{\partial}_{\varphi }^{\ast }+\bar{%
\partial}_{\varphi }^{\ast }\bar{\partial}:\Lambda ^{0,0}(C,\mathcal{O}%
(N))\rightarrow \Lambda ^{0,0}(C,\mathcal{O}(N))
\end{equation}

to be the Laplacian operator with respect to the metric $g_{\varphi }=\exp
\varphi g_{0}$. Since $\Delta _{\varphi }$ is a well-defined elliptic
operator, its regularized determinant can be defined by the method of zeta
regularization. The following functional gives the difference of the
determinant with respect to the conformal change of the metric:

\begin{definition}
\begin{equation}
A_{n}(\varphi )=\log (\frac{\det \Delta _{\varphi }}{\det \Delta _{0}}).
\end{equation}
\end{definition}

In order to present $A_{n}$ more concretely, we consider the kernel of $%
\Delta $, which consists of all the polynomials of $z$ of degree at most $n$%
; thus; it is a complex linear space of dimension $n+1$ and has an
othonormal basis with respect to the $g_{0}$ metric as follows:
\begin{equation}
\alpha _{i}=\sqrt{(n+1)\binom{n}{i}}z^{i},\text{ \ \ \ \ \ }i=1,\cdots ,n.
\end{equation}

Because $C$ is one dimensional, by the Serre Duality, $A\left( \varphi
\right) $ can also be viewed as analytic torsion associated to chain complex
(\ref{chain}).\ Thus, by the the anomaly formula of Bismut-Gillet-Soul\'{e}
\cite{BGS1}, the following variational formula is obtained by Gillet-Soul%
\'{e} \cite{GiSo1}:

\begin{proposition}
\begin{equation}
A_{n}(\varphi )=-{\frac{1}{{2}}}\int {|\nabla \varphi |^{2}\mu }-(n+1)\int {%
\varphi \mu }+\log \det {(<\alpha _{i},\alpha _{j}>}_{\varphi }{)_{n+1,n+1}.}
\end{equation}
\end{proposition}

Notice the similarity of $A_{n}$ with the left hand side of the
Moser-Trudinger Inequality (\ref{Moser-Trudinger}), which can be viewed as
the log determinant of the scalar Laplacian on the sphere (without the
involvement of the line bundle), and whose proof is related to the Yamabe
problem of $S^{2}.$

It is a simple observation that the value of this functional is invariant
when a constant function is added to $\varphi ;$ that is, for any constant $%
c\in \mathbb{R},$%
\begin{equation*}
A(\varphi )=A(\varphi +c).
\end{equation*}%
Hence, by subtracting a suitable constant, we may assume that $\int \varphi
\mu =0.$

For future convenience, we define the following functional%
\begin{equation}
B_{n}(\varphi )=\log \det {(<\alpha _{i},\alpha _{j}>}_{\varphi }{%
)_{n+1,n+1},}
\end{equation}%
which is the L$^{2}$ contribution of the anomaly formula; this is also the
fully non-linear term of the functional $A_{n}$.

In \cite{GiSo1}, by direct estimate following Moser, Gillet and Soul\'{e}
proved \ the following

\begin{theorem}
(\textbf{Gillet-Soul\'{e}})\label{oldcase}If $\varphi $ is rotationally
symmetric (with respect to a given coordinate system on $C=S^{2}$), then
\begin{equation*}
B_{n}(\varphi )\leq (\frac{1}{2}-\varepsilon _{n})\int {|\nabla \varphi
|^{2}\mu }+(n+1)\int {\varphi \mu +C}_{n}{,}
\end{equation*}%
where $\varepsilon _{n}>0$ and $C_{n}$ are constants depending only on $n;$
in particular, $A(\varphi )$ is bounded by a constant independent of $%
\varphi $.
\end{theorem}

\section{Bounded-ness of the functional $A(\protect\varphi )$%
\setcounter{equation}{0} \setcounter{theorem}{0}}

The main goal of this Section is to prove Theorem \ref{main}, which is a
generalization of Theorem~\ref{oldcase}. Our strategy is to simplify the
general case to the case which has been treated by Gillet and Soul\'{e}.\
Our approach involves symmetrization of various geometric quantities and
their estimates, where the difficulty is to treat the highly non-linear term
$B_{n}.$

First, we prove the following algebraic lemma.

\begin{lemma}
\label{algebra}Let $d\lambda $ be a finite positive measure of space $C$.
Let $0\leq a_{0}<\cdots \leq a_{n}$ be $n+1$ real numbers. For any $\sigma
\in S_{n+1}$, define
\begin{equation}
S_{\sigma ,\lambda }=S_{\sigma }=\prod_{i=0}^{n}\int_{C}|z|^{a_{i}+a_{\sigma
(i)}}d\lambda ,
\end{equation}%
then we have
\begin{equation}
S_{\sigma }\leq S_{id}.
\end{equation}
\end{lemma}

\begin{proof}
This is just an extension of the Cauchy-Schwartz Inequality. We run an
induction argument on $n$. The statement is obviously true for $n=1,$ in
which case it is just the usual Cauchy-Schwartz inequality. Suppose it is
true for $1,\cdots ,n-1$. Fix a $\sigma \in S_{n}$. If $\sigma $ can be
written as product of two shorter cycles, we can apply the induction
hypothesis. Without loss of generality, we may assume
\begin{equation*}
\sigma =(0,i_{1},i_{2},\cdots ,i_{n}).
\end{equation*}%
To further simplify the notation, we assume that $a_{i}=i.$\ (Notice that
the proof of general case follows similarly). Let
\begin{equation*}
\sigma ^{\prime }=(i_{1},i_{2},\cdots ,i_{n})
\end{equation*}%
be a reduced $(n-1)$-cycle. Then we have,
\begin{equation}
\log ({S_{\sigma }/S_{\sigma ^{\prime }}})=\log \int {|z|^{i_{1}}d\lambda }%
+\log \int {|z|^{i_{n}}d\lambda }-\log \int {|z|^{i_{1}+i_{n}}d\lambda }%
-\log \int {1\ d\lambda .}
\end{equation}%
Simple calculation shows that
\begin{equation}
u(t)\overset{\text{def}}{=}\log \int {|z|^{i_{1}+i_{n}-t}d\lambda }+\log
\int {|z|^{t}d\lambda }
\end{equation}%
is an increasing function for $t\in \lbrack 0,\frac{i_{0}+i_{n}}{2}]$. Apply
this fact to $0=t_{1}<t_{2}=min\{i_{1},i_{n}\}\leq \frac{i_{1}+i_{n}}{2}$,
we get
\begin{equation*}
S_{\sigma }\leq S_{\sigma ^{\prime }}.
\end{equation*}%
Now we can apply the induction hypothesis to prove that
\begin{equation*}
S_{\sigma }\leq S_{id}.
\end{equation*}%
Thus we have finished the induction proof.
\end{proof}

Lemma~\ref{algebra} allows us to estimate the terms appearing in the
determinant $B_{n}$ by the diagonal term. Now we state the main theorem of
this section:

\begin{theorem}
\label{main-2}%
\begin{equation}
B_{n}(\varphi )\leq (\frac{1}{2}-\varepsilon _{n})\int {|\nabla \varphi
|^{2}\mu }+(n+1)\int {\varphi \mu +C}_{n}{;}
\end{equation}%
where $\varepsilon _{n}>0$ and $C_{n}$ are constants depending only on $n;$
in particular, $A(\varphi )$ is bounded by a constant independent of $%
\varphi $.
\end{theorem}

Theorem~\ref{main} is then an easy corollary of Theorem~\ref{main-2}.

From an analytic point of view, Theorem \ref{main-2} can be viewed as a
\textit{multi-particle Moser-Trudinger}\textbf{\ }\textit{Inequality},
generalizing the original inequality (\ref{Moser-Trudinger}) of Moser and\
Trudinger~\cite{Moser-1}. We comment that the cofficient of the energy in (%
\ref{unknown}) has a uniform upper bound independent of $n$, making it a
very precise inequality. In the next sections, we will also discuss its
sharp form.

Finally, we give the prove of Theorem \ref{main-2}.

\begin{proof}
First, by Lemma~\ref{algebra} and the definition of the determinant, we have
\begin{equation}
|\det (<\alpha _{i},\alpha _{j}>_{\varphi })|\leq \sum_{\sigma \in
S_{n+1}}S_{\sigma ,\lambda }\leq n!S_{id,\lambda },  \label{jay1}
\end{equation}%
where $d\lambda ={\frac{\exp (\varphi )}{{(1+|z|^{2})^{n}}}}\mu $.

Notice that $l=\{z:|z|=1\}$ is a closed geodesic on $C,$ with respect to the
standard metric. Because of the rotational invariance of $C=S^{2},$ with out
loss of generality, we assume that $\max_{z\in C}\varphi (z)$ is obtained at
a point on $l.$ Then, we discuss the rotational re-arrangement of $\varphi ,$
following \cite{SY}. Define $\varphi ^{\ast }$ as the rotationally symmetric
increasing re-arrangement of $\varphi $ for $|z|\leq 1$ and the rotationally
symmetric decreasing re-arrangement function of $\varphi $ for $|z|\geq 1.\ $%
Hence, we have%
\begin{eqnarray*}
\varphi ^{\ast }(0) &=&\min_{\left\vert z\right\vert \leq 1}\varphi (z), \\
\varphi ^{\ast }(\infty ) &=&\min_{\left\vert z\right\vert \geq 1}\varphi
(z), \\
\varphi (1) &=&\max_{S^{2}}\varphi (z).
\end{eqnarray*}%
It is easy to see that $\varphi ^{\ast }$ is continuous. By \cite{SY},%
\begin{eqnarray}
\int_{C}|\nabla \varphi |^{2}\mu &\geq &\int_{C}|\nabla \varphi ^{\ast
}|^{2}\mu ,  \label{jay4} \\
\int \varphi \mu &=&\int_{C}\varphi \mu .  \label{jay5}
\end{eqnarray}%
Hence, $\varphi ^{\ast }\in C^{0}(C)\cap W^{1,2}(C).$

For $\left\vert z\right\vert \geq 1,$ define, for $i=0,\cdots ,n,$%
\begin{eqnarray*}
P_{i} &=&\frac{|z|^{2i}}{(1+\left\vert z\right\vert ^{2})^{n}}, \\
P_{i}^{\prime } &=&\frac{1}{(1+\left\vert z\right\vert ^{2})^{n}}.
\end{eqnarray*}%
Then, for $\left\vert z\right\vert \geq 1,$ we have%
\begin{equation}
P_{i}\leq P_{i}^{\prime }\leq 2^{2i}P_{i}\leq 2^{2n}P_{i}.
\end{equation}%
Notice that both $P_{i}^{\prime }$ and $\varphi ^{\ast }$ are non-increasing
for $\left\vert z\right\vert \geq 1,$ it is easy to check that%
\begin{eqnarray}
\log \int_{\left\vert z\right\vert \geq 1}P_{i}\exp (\varphi )\mu &\leq
&\log \int_{\left\vert z\right\vert \geq 1}P_{i}^{\prime }\exp (\varphi )\mu
\leq \log \int_{\left\vert z\right\vert \geq 1}P_{i}^{\prime }\exp (\varphi
^{\ast })\mu  \notag \\
&\leq &\log \int_{\left\vert z\right\vert \geq 1}P_{i}\exp (\varphi ^{\ast
})\mu +2n\log 2.  \label{jay2}
\end{eqnarray}%
By the reflection symmetry of $C=S^{2}$ and the definition of $\varphi
^{\ast },$ similarly we have%
\begin{equation}
\log \int_{\left\vert z\right\vert \leq 1}P_{i}\exp (\varphi )\mu \leq \log
\int_{\left\vert z\right\vert \leq 1}P_{i}\exp (\varphi ^{\ast })\mu +2n\log
2.  \label{jay3}
\end{equation}%
By (\ref{jay1}) (\ref{jay2}) and (\ref{jay3}), we have
\begin{equation}
\log \det {(<\alpha _{i},\alpha _{j}>_{\varphi })_{n+1,n+1}}\leq n!\log \det
{(<\alpha _{i},\alpha _{j}>_{\bar{\varphi}})_{n+1,n+1}}.  \label{jay10}
\end{equation}%
Finally, combine (\ref{jay4}), (\ref{jay5}) and (\ref{jay10}), we have proved%
\begin{equation}
A(\varphi )\leq A(\bar{\varphi})+C(n).
\end{equation}%
Finally, by applying Theorem~\ref{oldcase}, we have proved the theorem.
\end{proof}

\section{$0$ is a local maximum}

\label{sec:3}\setcounter{equation}{0} \setcounter{theorem}{0}Now that the
functional $A_{n}$ is bounded by Theorem \ref{main-2}, it is natural to
examine the sharp upper-bound of $A_{n}.$ In this Section we show that the
functional $A_{n},$ has a local maximal point at constant function, $0,$ for
any fixed $n.$\ This stability result indicates that the maximum of
functional $A_{n}$ should be achieved at the constant function.

We begin with some notations. For simplicity, we define the following local $%
(n+1)\times (n+1)$ matrices: for $\varphi ,f\in C^{\infty }(C,\mathbb{R}),$
write%
\begin{eqnarray*}
(M_{\varphi })_{ij} &=&<\alpha _{i},\alpha _{j}>_{\varphi }, \\
M_{\varphi }^{\prime }(f) &=&<f\alpha _{i},\alpha _{j}>_{\varphi }, \\
M_{\varphi }^{\prime \prime }(f) &=&<f^{2}\alpha _{i},\alpha _{j}>_{\varphi
}, \\
(M_{\varphi })^{ij} &=&(M_{\varphi })^{-1}.
\end{eqnarray*}

First, we give the following first variation formula:

\begin{lemma}[First Variation of $A_{n}$]
If $g_{\varphi }$ is an extreme metric for the functional $A_{n},$ then%
\begin{equation}
\limfunc{tr}(M_{\varphi }^{-1}[\frac{\alpha _{i}\overline{\alpha _{j}}}{{(1+}%
\left\vert z\right\vert {^{2})^{n}}}{]}_{ij})+\Delta \varphi -(n+1)=0.
\label{Euler}
\end{equation}%
In particular, the constant function $0$ is a critical point of $A_{n}$ for
any $n.$
\end{lemma}

\begin{proof}
The derivation of the above Euler-Lagrange equation is standard, which we
omit here.\ To verify that $0$ is critical, use the fact that $M_{0}=Id,$ by
(\ref{Euler}), we only need to show%
\begin{equation}
n+1=\sum_{i=0}^{n}\frac{\alpha _{i}\cdot \overline{\alpha _{i}}}{%
(1+\left\vert z\right\vert ^{2})^{n}};
\end{equation}%
or,%
\begin{equation*}
n+1=\sum_{i=0}^{n}\frac{\left( n+1\right) \binom{n}{i}\left\vert
z\right\vert ^{2i}}{(1+\left\vert z\right\vert ^{2})^{n}},
\end{equation*}%
which is obvious.
\end{proof}

We then proceed to compute the second variation of $A_{n}.$

\begin{proposition}
For $\varphi ,$ $f\in C^{\infty }(C,\mathbb{R}),$ the following holds%
\begin{equation}
\frac{d^{2}}{d^{2}t}|_{t=0}A(\varphi +tf)=\limfunc{tr}[M_{\varphi
}^{-1}M_{\varphi }^{\prime \prime }(f)]-\limfunc{tr}[M_{\varphi
}^{-1}M_{\varphi }^{\prime }\left( f\right) M_{\varphi }^{-1}M_{\varphi
}^{\prime }\left( f\right) ]-\int \left\vert \nabla f\right\vert ^{2}\exp
\varphi \mu .  \label{key-3}
\end{equation}
\end{proposition}

The proof of (\ref{key-3}) is standard, so we omit it here.

Finally, we state the main result of this section.

\begin{theorem}
\label{0localmax}$0$ is a local maximal point of the functional $A_{n}.$
\end{theorem}

\begin{proof}
By (\ref{key-3}), we need to show the following inequality for any $f\in
C^{\infty }(C),$%
\begin{equation}
\sum_{i}\int_{C}f^{2}\alpha _{ij}\mu -\sum_{i,j}|\int_{C}f\alpha _{ij}\mu
|^{2}\leq \int \left\vert \nabla f\right\vert ^{2}\mu ,  \label{key-1}
\end{equation}%
where, for notational simplicity, we denote%
\begin{equation*}
\alpha _{ij}=\frac{\alpha _{i}\overline{\alpha _{j}}}{{(1+}\left\vert
z\right\vert {^{2})^{n}}}{.}
\end{equation*}

Using the fact that $\left\{ \alpha _{i}\right\} $ is an orthonormal basis,
the left hand side of (\ref{key-1}) can be re-written as%
\begin{eqnarray}
&&\dint\limits_{x\in C}\dint\limits_{y\in C}\sum_{ij}(f^{2}(x)\alpha _{ij}(x)%
\overline{\alpha _{ij}(y)}-f(x)f(y)\alpha _{ij}(x)\overline{\alpha _{ij}(y)}%
)\mu _{y}\mu _{x}  \notag \\
&=&\diint\limits_{C\times C}\frac{1}{2}\sum_{i,j}\alpha _{ij}(x)\overline{%
\alpha _{ij}(y)}(f(x)-f(y))^{2}\mu _{x}\mu _{y}  \notag \\
&=&\diint\limits_{C\times C}\frac{(n+1)^{2}}{2}(f(x)-f(y))^{2}\sum_{i,j}%
\binom{n}{i}\binom{n}{j}\frac{z_{x}^{i}\overline{z}_{x}^{j}z_{y}^{i}%
\overline{z}_{y}^{j}}{(1+|z_{x}|^{2})^{n}(1+|z_{y}|^{2})^{n}}\mu _{x}\mu _{y}
\notag \\
&=&\diint\limits_{C\times C}\frac{(n+1)^{2}}{2}(f(x)-f(y))^{2}(\frac{%
\left\vert 1+z_{x}\overline{z_{y}}\right\vert ^{2}}{%
(1+|z_{x}|^{2})(1+|z_{y}|^{2})})^{n}\mu _{x}\mu _{y}  \notag \\
&=&\diint\limits_{C\times C}\frac{(n+1)^{2}}{2}(f(x)-f(y))^{2}\cos ^{2n}(%
\sqrt{\pi }d_{\mu }(x,y))\mu _{x}\mu _{y},  \label{key-4}
\end{eqnarray}%
where $d_{\mu }(\cdot ,\cdot )$ is the distance function of $C=S^{2}$ with
respect to the metric $\mu .$ To do the local computation, sometimes we
switch to the standard metric of $S^{2}\subset \mathbb{R}^{3}$, which has
radius one and volume form $dv=4\pi \mu .$ If we use $d(\cdot ,\cdot )$ to
denote the distance function with respect to the standard metric, then
\begin{equation}
d=\sqrt{4\pi }d_{\mu }.
\end{equation}%
We continue our computation. By (\ref{key-4}), (\ref{key-1}) is equivalent
to the following%
\begin{equation}
\diint\limits_{C\times C}\frac{(n+1)^{2}}{2}(f(x)-f(y))^{2}\cos ^{2n}(\frac{%
d(x,y)}{2})\mu _{x}\mu _{y}\leq \int \left\vert \nabla f\right\vert ^{2}\mu .
\label{key-5}
\end{equation}%
Let $dl$ be the induced measure on embedded curves in $C$ from $dv.$\ For
fixed points $x,y\in C,$ let $l(x,y)$ be the shortest geodesic connecting $x$
and $y.$ Then, the following estimate holds%
\begin{eqnarray}
\left\vert f(x)-f(y)\right\vert ^{2} &\leq &(\dint\limits_{l(x,y)}|\nabla
f(p)|\cdot |\cos \delta (p,x,y)|dl_{p})^{2}  \notag \\
&\leq &d(x,y)\dint\limits_{l(x,y)}|\nabla f(p)|^{2}\cos ^{2}\delta
(p,x,y)dl_{p},  \label{key-6}
\end{eqnarray}%
where $\delta (p,x,y)$ denotes the angle between the direction of $l(x,y)$
and $\nabla f(p).$\ Obviously, $\delta $ can be determined only by $p$ and $%
x $; hence, it can also be written as $\delta (p,x)$.

Fix $x\in C$ and write $\mu _{y}=\frac{1}{4\pi }\sin \theta d\theta d\varphi
,$ with $\theta =d(x,y),$ $\varphi \in \lbrack 0,2\pi ),\theta \in \lbrack
0,\pi ]$ being the spherical coordinates.\ Hence, if $l=d(p,x),$%
\begin{eqnarray}
&&\diint\limits_{C\times C}|f(x)-f(y)|^{2}\cos ^{2n}(\frac{d(x,y)}{2})\mu
_{x}\mu _{y}  \notag \\
&\leq &\int_{x}\frac{1}{4\pi }\int_{0}^{2\pi }\int_{0}^{\pi }\cos ^{2n}(%
\frac{\theta }{2})\theta \sin \theta \int_{0}^{\theta }|\nabla f(p)|^{2}\cos
^{2}\delta (p,x)dl_{p}d\theta d\varphi \mu _{x}  \label{key-7} \\
&=&\frac{1}{4\pi }\int_{x}\int_{0}^{2\pi }\int_{0}^{\pi }(\int_{l}^{\pi
}\cos ^{2n}(\frac{\theta }{2})\theta \sin \theta d\theta )\cos ^{2}\delta
(p,x)\left\vert \nabla f(p)\right\vert ^{2}dld\varphi \mu _{x}.  \notag
\end{eqnarray}%
Thus, if we define%
\begin{equation}
G(l)=\int_{l}^{\pi }\cos ^{2n}(\frac{\theta }{2})\theta \sin \theta d\theta
),
\end{equation}%
The right hand side of (\ref{key-7}) can be written as%
\begin{eqnarray}
&&\frac{1}{4\pi }\int_{x}\int_{0}^{2\pi }\int_{0}^{\pi }\cos ^{2}\delta
(x,p)\left\vert \nabla f(p)\right\vert ^{2}\frac{G(l)}{\sin l}\sin l\text{ }%
dld\varphi \mu _{x}  \notag \\
&=&\frac{1}{4\pi }\int_{x}\int_{p}\int_{0}^{\pi }\cos ^{2}\delta
(x,p)\left\vert \nabla f(p)\right\vert ^{2}\frac{G(l)}{\sin l}\mu _{p}\mu
_{y} \\
&=&\int_{p}\left\vert \nabla f(p)\right\vert ^{2}\mu _{p}(\int \frac{%
G(d(x,p))}{\sin d(x,p)}\cos ^{2}\delta (x,p)\mu _{x}).
\end{eqnarray}%
Hence, to prove (\ref{key-1}), it is sufficient to prove the following%
\begin{equation}
\frac{\left( n+1\right) ^{2}}{2}\int \frac{G(d(x,p))}{\sin d(x,p)}\cos
^{2}\delta (x,p)\mu _{x}\leq 1.  \label{key-8}
\end{equation}

We parameterize $x\in C$ be $l=d(x,p)$ and $\delta =\delta (x,p),$ such that
$\mu _{x}=\frac{1}{4\pi }\sin l$ $dld\delta .\ $Hence, the left hand side of
(\ref{key-8}) equals

\begin{eqnarray*}
&&\frac{\left( n+1\right) ^{2}}{2}(\frac{1}{4}\int_{0}^{\pi }\frac{G(l)}{%
\sin l}\sin ldl) \\
&=&\frac{\left( n+1\right) ^{2}}{8}\int_{0}^{\pi }G(l)dl \\
&=&\frac{\left( n+1\right) ^{2}}{8}\int_{0}^{\pi }\cos ^{2n}(\frac{\theta }{2%
})\theta ^{2}d\theta .
\end{eqnarray*}%
By a change of variable $t=\frac{\theta }{2}$ and use the integration by
part trick, we get
\begin{equation*}
\frac{\left( n+1\right) ^{2}}{8}\int_{0}^{\pi }\cos ^{2n}(\frac{\theta }{2}%
)\theta ^{2}d\theta =2(n+1)\int_{0}^{\frac{\pi }{2}}t\cos ^{2n+2}t\text{ }dt.
\end{equation*}%
Define
\begin{equation}
J_{n}=2(n+1)\int_{0}^{\frac{\pi }{2}}t\cos ^{2n+2}t\text{ }dt.  \label{J}
\end{equation}%
Then (\ref{key-8}) is equivalent to the following claim:

\begin{claim}
\begin{equation*}
J_{n}<1,\text{ \ \ for }n\in \mathbb{N}\text{.}
\end{equation*}
\end{claim}

\begin{proof}
Using integration by part, it is easy to get the following recursive formula%
\begin{equation}
J_{n}=-\frac{1}{2n+2}+\frac{2n+1}{2n}J_{n-1}.
\end{equation}%
To prove the claim, we make the following statement%
\begin{equation}
J_{n-1}\geq \frac{3n+2}{3n+3}.  \tag{*}  \label{*}
\end{equation}%
Obviously (\ref{*}) is wrong for $n=1.$ If this is the case for all $n$,
then the claim is proved.\ Otherwise, let $N$ be the smallest positive
integer such that (\ref{*}) is true. Then by (\ref{key-9}),%
\begin{eqnarray*}
J_{N}-\frac{3N+5}{3N+6} &\geq &-\frac{1}{2N+2}+(\frac{2N+1}{2N})\frac{3N+2}{%
3N+3}-\frac{3N+5}{3N+6} \\
&=&\frac{4}{6N(N+1)(N+2)}>0.
\end{eqnarray*}%
Thus, easy to show by induction that (\ref{*}) holds for all $n\geq N.\ $%
Therefore, by (\ref{key-9}),%
\begin{eqnarray*}
J_{n}-J_{n-1} &=&\frac{J_{n-1}}{2n}-\frac{1}{2n+2}\geq \frac{3n+2}{(2n)(3n+3)%
}-\frac{1}{2n+2} \\
&=&\frac{1}{3n(n+1)}>0;
\end{eqnarray*}%
in other words, by induction, $J_{n}$ is monotone increasing for $n\geq N.$
On the other hand, applying a standard stationary phase argument to (\ref{J}%
), it is easy to see that
\begin{equation*}
\lim_{n\rightarrow \infty }J_{n}=1.
\end{equation*}%
Hence, for $n<N,$ $J_{n}\leq \frac{3n+2}{3n+3}<1;$ for $n\geq N,$ $J_{n}\leq
J_{n+1}\leq \cdots <1.$ The claim is proved.
\end{proof}

In summary, we have finished the proof of Theorem~\ref{0localmax}. It is
also clear from (\ref{key-6}) that the equality holds if and only if $%
\left\vert \nabla f\right\vert =0$ everywhere; \textit{i.e.}, $f$ is a
constant function.
\end{proof}

\section{Conjectured sharp inequality and its different formulations}

\setcounter{equation}{0} \setcounter{theorem}{0}Inspired by early
computations, we pose the sharp inequality of our problem, and point out its
links to different classical results.

\subsection{Conjecture on the sharp 2-D Szeg\"{o} limit theorem}

We now discuss the classical Szeg\"{o} limit theorem. Let $z$ be the complex
coordinate function on $\mathbb{C}$, $S^{1}\subset \mathbb{C}$ be the
standard circle of radius $1,$ and $ds=\frac{1}{2\pi }d|z|$ be the regular
measure on $S^{1}.$ For a function $\varphi \in C^{\infty }(S^{1},\mathbb{R)}
$, define
\begin{equation*}
\alpha _{ij}\overset{\text{def}}{=}z^{i}\overline{z}^{j}=z^{i-j},
\end{equation*}%
\begin{equation}
B_{n}(\varphi ,S^{1})\overset{\text{def}}{=}\log \det [\int_{S^{1}}\alpha
_{ij}(\exp \varphi )ds]_{(n+1)\times (n+1)}.
\end{equation}%
Szeg\"{o} proved the following

\begin{theorem}[Szeg\"{o} Limit Theorem]
\label{Szego}$B_{n}(\varphi ,S^{1})$ is non-decreasing with respect to $n.$
Furthermore,%
\begin{equation*}
B_{n}(\varphi ,S^{1})\leq {\frac{1}{{2}}}\int {|\nabla \varphi |^{2}ds}%
+(n+1)\int {\varphi ds.}
\end{equation*}%
Equality holds iff $\varphi $ is a constant.
\end{theorem}

Noticing the similarity, we pose the following conjecture for the two sphere:

\begin{conjecture}
\label{conjecture-2}For any $\varphi \in L^{2}(S^{2},\mathbb{R}),$ we have%
\begin{equation*}
A_{n}(\varphi )\leq 0,
\end{equation*}%
or, equivalently,%
\begin{equation}
B_{n}(\varphi )\leq {\frac{1}{{2}}}\int {|\nabla \varphi |^{2}\mu }%
+(n+1)\int {\varphi \mu ;}  \label{unknown}
\end{equation}%
equality holds iff $\varphi $ is a constant function.
\end{conjecture}

Conjecture~\ref{conjecture-2} equivalent to Conjecture~\ref{conjecture-2.0}.

\begin{remark}
This is because the coefficient of the energy in (\ref{unkown}) is twice
that of (\ref{Moser-Trudinger}), the original Moser-Trudinger inequality.
Thus, it is easy to see that the above conjecture is true for $n=0,1.$
\end{remark}

\begin{remark}
By Theorem~\ref{main-2}, Conjecture~\ref{conjecture-2} holds for $\varphi $
with large energy for general $n$. It also holds for $\varphi $ with very
small energy by Theorem~\ref{0localmax}.
\end{remark}

\begin{remark}
Due to the non-linear nature of our problem, for fixed $\varphi ,$ it is not
clear if $B_{n}\left( \varphi \right) $ is monotone with respect to $n$.
Notice that monotonicity is the key ingredient in Szeg\"{o}'s approach of
Theorem~\ref{Szego}.
\end{remark}

\begin{remark}
Okikiolu extended a weak form of Szeg\"{o} limit theorem to $S^{2}$ and $%
S^{3}.\ $In \cite{ko}, she proved a similar asymptotic results as that of
Szeg\"{o} for the determinant of the Toeplitz operator $M_{\exp \varphi }$
(see the next subsection for more details) restricted to different levels of
spherical harmonics, and got . However, Conjecture \ref{conjecture-2} is a
holomorphic extension, in which we consider only the \textquotedblleft
holomorphic\textquotedblright\ part of the spherical harmonics.
\end{remark}

This conjecture can also be viewed as the sharp form of the multi-particle
Moser-Trudinger Inequality proved in Theorem~\ref{main-2}. Actually, the
proof of the original (sharp) Moser-Trudinger Inequality goes as follows:
First, prove that an analogue of Theorem \ref{main} holds; \textit{i.e., the
functional has a universal upper bound.} Then, get the sharp upper bound by
analysis the Euler equiation as a PDE. However, this strategy is hard to
follow to prove Conjecture~\ref{conjecture-2} since the Euler equation (\ref%
{Euler}) for our multi-particle problem is highly non-linear and is not even
a PDE.

\subsection{Toeplitz operator point of view}

The functional $B_{n}$ has appeared in classical analysis in another form.
Boutet de Monvel and Guillemin have considered the quantization problem for
the Toeplitz operator in the framework of the pseudo-differential operator
\cite{BdM-G}.

Let $F=\exp \varphi .$ Let the Toeplitz operator $M_{F}$ be the
multiplication operator with respect to $F.$ $M_{F}$ acting on $\Gamma (%
\mathcal{O}\left( n\right) )$ can be viewed as a quantization of $F:$%
\begin{equation*}
M_{F}=F\cdot :\Gamma (\mathcal{O}\left( n\right) )\rightarrow \Gamma (%
\mathcal{O}\left( n\right) ).
\end{equation*}%
Since $\{\alpha _{i}\}_{i=0}^{n}$ forms a basis of $\Gamma (\mathcal{O}%
\left( n\right) )$,%
\begin{equation*}
B_{n}(\varphi )=\log \det \int_{C}<M_{F}\alpha _{i},\alpha _{j}>\mu
\end{equation*}%
is the so-called Toeplitz determinant$.$ A special case of a theorem of
Boutet de Monvel-Guillemin and a result of Uribe can be stated as follows:

\begin{theorem}
\label{Toeplitz}

\begin{enumerate}
\item (Boutet de Monvel-Guillemin) For $n>>1,$ $B_{n}(\varphi )$ has the
following asymptotic expansion:%
\begin{equation*}
B_{n}(\varphi )=(n+1)\int_{C}\varphi \mu +\sum_{i=1}^{N}\frac{D_{i}(\varphi )%
}{\left( n+1\right) ^{i}}+\frac{\widetilde{D_{N+1}}(n,\varphi )}{\left(
n+1\right) ^{N+1}},
\end{equation*}%
where $N\in \mathbb{N}$, $D_{i}(\varphi )$'s are functionals independent of $%
n$ and dependent of $\varphi $ only, $\widetilde{D_{N+1}}(\varphi )$ is
bounded independent of $n$;

\item (Uribe, \cite{uribe-1}) furthermore,%
\begin{equation*}
D_{0}(\varphi )=\frac{1}{2}\int_{C}\left\vert \nabla \varphi \right\vert
^{2}\mu .
\end{equation*}
\end{enumerate}
\end{theorem}

Theorem~\ref{Toeplitz} shows that Conjecture~\ref{conjecture-2}, with a
fixed $\varphi $, is asymptotically correct for $n>>1$.

Notice that Uribe also gave an inductive method to compute all the $%
D_{i}(\varphi )^{\prime }$s. However, it is not clear that the complete
asymptotic expansion would lead to the sharp Moser Inequality.

\subsection{Probability formulation of functional $B_{n}$}

Using symmetries of $C=\mathbb{CP}^{1}=S^{2},$ we can reformulate the
problem and give it a probabilistic interpretation. More precisely, we prove
the following

\begin{proposition}
\label{prob}For $(z_{0},\cdots ,z_{n+1})\in (C)^{n+1},$ define%
\begin{equation}
K_{n}(z_{0},\cdots ,z_{n+1})=\frac{(n+1)^{n+1}\binom{n}{0}\cdots \binom{n}{n}%
}{(n+1)!}\dprod\limits_{0\leq i<j\leq n}(\sin \frac{d(z_{i},z_{j})}{2})^{2}.
\label{4-0}
\end{equation}%
Then,%
\begin{eqnarray}
\dint\limits_{C^{n+1}}K_{n}\prod \mu _{z_{i}} &=&1,  \label{4-1} \\
B(\varphi ) &=&\log \dint\limits_{C^{n+1}}\exp (\varphi (z_{0})+\cdots
\varphi (z_{n}))K_{n}\prod \mu _{z_{i}}.  \label{4-2}
\end{eqnarray}
\end{proposition}

\begin{proof}
It is easy to see that (\ref{4-1}) follows from (\ref{4-2}).\ Thus, we only
need to prove the latter. For simplicity, we write $\mu _{i}=\mu _{z_{i}}.$%
\begin{eqnarray*}
\exp B(\varphi ) &=&\det \left[ \int_{C}\alpha _{ij}\exp \varphi \mu \right]
\\
&=&\frac{(n+1)^{n+1}\binom{n}{0}\cdots \binom{n}{n}}{(n+1)!}%
\dint\limits_{C^{n+1}}\exp (\sum \varphi (z_{i}))\frac{\dsum\limits_{\sigma
\in S_{n}}(-1)^{\left\vert \sigma \right\vert }\prod z_{i}^{i}\overline{z}%
_{i}^{\sigma (i)}}{\prod (1+\left\vert z_{i}\right\vert ^{2})^{n}}\prod \mu
_{i} \\
&=&\frac{(n+1)^{n+1}\binom{n}{0}\cdots \binom{n}{n}}{(n+1)!}%
\dint\limits_{C^{n+1}}\exp (\sum \varphi (z_{i}))\frac{\dprod\limits_{i<j}%
\left\vert z_{i}-z_{j}\right\vert ^{2}}{\prod (1+\left\vert z_{i}\right\vert
^{2})^{n}}\prod \mu _{i}.
\end{eqnarray*}%
(\ref{4-2}) then will be obvious considering the following simple fact on $%
S^{2}:$%
\begin{equation}
\frac{\left\vert z_{i}-z_{j}\right\vert ^{2}}{(1+\left\vert z_{i}\right\vert
^{2})(1+\left\vert z_{j}\right\vert ^{2})}=\sin ^{2}\frac{d(z_{i},z_{j})}{2}.
\end{equation}%
We have finished the proof.
\end{proof}

From Proposition \ref{prob}, our study of the functionals $A$ and $B$ can be
viewed probabilistically. Considering the $(n+1)$-random variables $%
\{z_{i}\}_{i=0}^{n}$ on the standard sphere $S^{2},$ with joint distribution
given by $K_{n},$%
\begin{equation}
\exp B_{n}\left( \varphi \right) =E\left( \exp \left( \sum \varphi
(z_{i}\right) \right) ,
\end{equation}%
is the expectation of the random variable $\exp \left( \sum \varphi
(z_{i}\right) ).$\ Hence, for $n$ very large, our Conjecture \ref%
{conjecture-2} can also be viewed as a \textit{limit theorem} of the given
probability problem. This direction is partly motivated by the works of Kac~%
\cite{Kac-2} and Johansson~\cite{Johansson-1} on the probabilistic approach to the Szeg\"{o} limit theorem.

\bibliographystyle{plain}
\bibliography{bib,NSF1}

\end{document}